\title{Almost-tiling the plane by ellipses\thanks{Research of K. Kuperberg
was supported by NSF grant No. DMS-9704958.
Research of W. Kuperberg
was supported by NSF grant No. DMS-9704319.
Research 
of J. Matou\v sek was supported by
ETH Z\"urich, by
Czech Republic Grant GA\v{C}R 0194/1996, and by Charles University grants
 No. 193,194/1996. Research of P. Valtr was supported by DIMACS Center, by
Czech Republic Grant GA\v{C}R 0194/1996, and by Charles University grants
 No. 193,194/1996.}}

\newcommand{\cmt}[1]{\ifhmode\newline\fi{\sf *** \ \ #1 \\}}

\documentstyle[11pt]{article}

\newtheorem{theorem}{Theorem}

\newcommand{\Rd}{{\R^d}}
\newcommand{\R}{{\hbox{\rm I}\!\hspace{-0.025em}\hbox{\rm R}}}

\newcommand{\eps}{\varepsilon}

\newcommand{\CC}{{\cal C}}
\newcommand{\EE}{{\cal E}}

\newcommand{\heading}[1]{\vspace{1ex}\par\noindent{\bf #1}}
\newcommand{\subh}[1]{\noindent{\em #1}}

\newcommand\epsval{{10^{-5}}}
\newcommand\alphaval{{\frac\pi{16}}}
\newcommand\betaval{{\frac1{16}}}
\newcommand\Ifrac{{\frac\beta{16}}}
\newcommand\newsqs{{\frac\beta{20}}}
\newcommand\smallrt{{\frac\beta{80}}}
\newcommand\bigrt{{\frac18}}
\newcommand\newarc{{\frac\alpha{8}}}
\newcommand\maxrprime{{\frac1{20}}}

\long\def\onefigure#1#2{
\begin{figure*}[htb]
\begin{center}
#1
\end{center}
\caption{#2}
\end{figure*}
} 

\newcommand{\lipefig}[2]  
{\onefigure{\def\IPEfile{packing-#1.ipe}\input{\IPEfile}}{\label{f:#1} #2} }

\author{
{\sc Krystyna Kuperberg}\\ 
  {\footnotesize  Department of Mathematics}\\[-1.5mm]
  {\footnotesize  Auburn University}\\[-1.5mm]
  {\footnotesize   Auburn, AL 36849-5310}\\[-1.5mm]
  {\footnotesize   U.S.A.}
\\[-1.5mm]   {\footnotesize e-mail: {\tt kuperkm@mail.auburn.edu}}
\and{\sc Wlodzimierz Kuperberg}\\
  {\footnotesize  Department of Mathematics}\\[-1.5mm]
  {\footnotesize  Auburn University}\\[-1.5mm]
  {\footnotesize   Auburn, AL 36849-5310}\\[-1.5mm]
  {\footnotesize   U.S.A.}
\\[-1.5mm]   {\footnotesize e-mail: {\tt kuperwl@mail.auburn.edu}}
\and{\sc Ji\v{r}\'{\i} Matou\v{s}ek}\\
   {\footnotesize Department of Applied Mathematics}\\[-1.5mm]
   {\footnotesize  Charles University}\\[-1.5mm]
   {\footnotesize  Malostransk\'{e} n\'{a}m. 25, 118~00~~Praha~1}\\[-1.5mm]
   {\footnotesize Czech Republic}\\[-1.5mm]
   {\footnotesize and}\\[-1.5mm]
   {\footnotesize Institut f\"ur Theoretische Informatik }\\[-1.5mm]
   {\footnotesize ETH Z\"urich, Switzerland }
\\[-1.5mm]   {\footnotesize e-mail: {\tt matousek@kam.mff.cuni.cz}}
   \and
   {\sc Pavel Valtr}\\
   {\footnotesize Department of Applied Mathematics}\\[-1.5mm]
   {\footnotesize  Charles University}\\[-1.5mm]
   {\footnotesize  Malostransk\'{e} n\'{a}m. 25, 118~00~~Praha~1}\\[-1.5mm]
   {\footnotesize Czech Republic}\\[-1.5mm]
   {\footnotesize and}\\[-1.5mm]
   {\footnotesize DIMACS Center }\\[-1.5mm]
   {\footnotesize Rutgers University, NJ, U.S.A.} 
 \\[-1.5mm]   {\footnotesize e-mail: {\tt valtr@kam.mff.cuni.cz}}    }
\date{}

\begin{document}

\maketitle

\centerline{\em Dedicated to our friend Imre B\'ar\'any on the occasion of his
50-th birthday.}
\begin{abstract}
\noindent
For any  $\lambda>1$ we construct a periodic and locally finite
packing of the plane with ellipses whose
$\lambda$-enlargement covers the whole plane.
This answers a question of Imre B\'ar\'any.
On the other hand, we show that if $\CC$ is a packing in the plane
with circular discs of radius at most $1$, then its $(1+\epsval)$-enlargement
covers no square with side length $4$.

\end{abstract}

\section{Introduction}
Let $\CC$ be a system (finite or infinite) of centrally symmetric
convex bodies in $\Rd$ with disjoint
interiors; we call  such a $\CC$ a {\em packing\/}.
For a real number $\eps>0$ and for $C\in \CC$, we let
$C^\eps$ denote $C$ enlarged by
the factor $1+\eps$ from its center, that is,
$C^\eps =(1+\eps)(C-c)+c$, where $c$ stands for the center of
symmetry $C$. Let us call the closure of the
set $C^\eps\setminus C$ the {\em $\eps$-ring\/} of $C$.
We call the system $\CC^\eps=\{C^\eps;\, C\in\CC\}$ 
the {\em $(1+\eps)$-enlargement} of $\CC$.

For a class $\CC_0$ of centrally symmetric convex bodies in
$\Rd$, we define the {\em B\'ar\'any number\/} of $\CC_0$ as
the infimum of the numbers $\eps>0$ such that there exists
a packing $\CC\subseteq \CC_0$ whose $(1+\eps)$-enlargement
covers the whole $\Rd$.

We observe that the system of all circular discs in the plane
has B\'ar\'any number 0, since we can produce the desired
packing for any $\eps>0$ by using larger and larger discs
(add discs to the packing one by one,
and in the $i$th step, choose the $i$th disc so that its
$\eps$-ring covers the disc of radius $i$ around the origin).
A different situation may occur if the diameter of the
bodies in $\CC_0$ is bounded.

Motivated by a problem concerning convex polytopes, Imre B\'ar\'any~\cite{imre}
raised a problem which in our terminology can be rephrased
as follows: If $\EE$ stands for the class of all ellipses
of diameter at most 1, is the B\'ar\'any
number of $\EE$ zero?

In this paper we give a positive answer to this question:

\begin{theorem}\label{t:ell}
For every $\lambda>1$ there is a periodic
packing of the plane with ellipses whose
$\lambda$-enlargement is a covering.
\end{theorem}

On the other hand, if we allow only discs
of bounded radius, then B\'ar\'any's question has a negative answer:

\begin{theorem}\label{t:discs}
Let $\CC$ be a packing of the plane with circular discs of radius
at most $1$. Then $(1+\epsval)$-enlargement 
of $\CC$ covers no square with side length $4$.
\end{theorem}

\heading{Remarks. }
Our packing in Theorem~\ref{t:ell} is locally finite and 
the details in the construction can be done so that all ellipses
in the packing have diameter between $\eps/10$ and $1$,
where $\eps=\lambda-1$
(however, their width varies from 
$\exp [-c/\varepsilon\log^2(1/\varepsilon )]$ to $c$ 
and we need $\exp [c/\varepsilon\log(1/\varepsilon )]$ of them
on each unit square).
Our methods can be used to prove that Theorem~\ref{t:discs}
(possibly with a different positive constant instead of $\epsval$)
holds also in any dimension $d\ge 2$ and when
$\CC_0$ consists of convex bodies in $\Rd$ with a constant-bounded 
diameter and curvature.
We do not prove these generalizations here, since the
idea remains the same but the details become messy.
The value $\eps=\epsval$ in Theorem~\ref{t:discs}
is certainly not the best possible one could get by
our proof method, but it seems that a different method would be
needed to determine the B\'ar\'any number for discs in the plane exactly or at
least to prove a reasonable lower bound for its value.

\section{Almost-tiling by ellipses}
Throughout the construction, a number $\lambda >1$ remains
fixed. Choose an integer $n$ such that the regular
$2n$-polygon $P=P_{2n}$ circumscribed about a circular disk
$D$ is contained in the interior of the $\lambda$-enlargement
of $D$. Denote two antipodal vertices of $P$ by $v^-$
and $v^+$. Suppose $T$ is a triangle with a horizontal base
$B$ and vertex $v$ above $B$. Then there is a (unique)
polygon $P(T)$ satisfying the following conditions (see
Fig. 1 which illustrates the case $n=4$):

(i) $P(T)\subset T$;
(ii) There is an affine transformation $A$ such
that $A(P)=P(T)$;
(iii) $A(v^-)$ is the midpoint of $B$ and $A(v^+)=v$;
(iv) The angles at $v$ of $P(T)$ and of $T$ are equal.

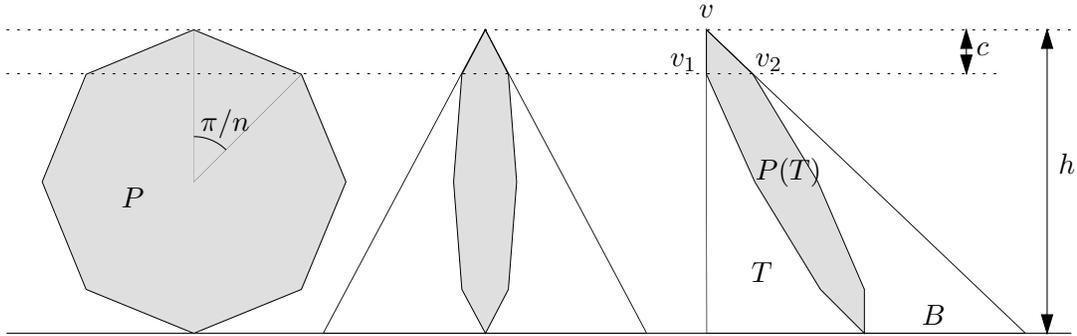
\lipefig{1}{The affine copy of $P$ properly inscribed in $T$.}

We say that $P(T)$ is an {\it affine copy of $P$ properly
inscribed in $T$.} The existence and uniqueness of $P(T)$
follows from the fact that $P$ itself, oriented so that one
of its main diagonals is vertical, is properly inscribed in
a triangle, and an affine transformation that sends this
triangle onto $T$ (top vertex onto top vertex and base onto
base) determines $P(T)$ uniquely. Obviously, $P(T)$ contains
an inscribed ellipse, namely $A(D)$, whose
$\lambda$-enlargement contains a neighborhood of $P(T)$.

Observe the following property of the polygon $P(T)$:

(1) Let $v_1$ and $v_2$ be the vertices of $P(T)$
adjacent to $v$. Then the line $v_1v_2$ is parallel to $B$
and partitions the height $h$ of $T$ at the ratio of 
$$c:(h-c)=
\left(1-\cos{\pi\over n}\right):
\left(1+\cos{\pi\over n}\right),$$
where $c$ is the portion of $h$ containing $v$ (see Fig. 1).

It follows that

(2) For every vertex $w\neq v$ of $P(T)$ the
distance from $w$ to the line of $B$ is less than or equal to
$\mu h$, where $\mu<1$ is a positive constant independent from
$T$. Specifically,
$$\mu={1\over 2}\left(1+\cos{\pi\over n}\right).$$ 

The construction continues with the following lemma:
\medskip

{\bf Lemma. } {\em If $U$ is a neighborhood of a side of a
triangle $T$, then there is a polygonal region $W$ containing
$T\setminus U$ and contained in $T$, which can be tiled
by a finite collection of affine copies of $P$.}
\medskip
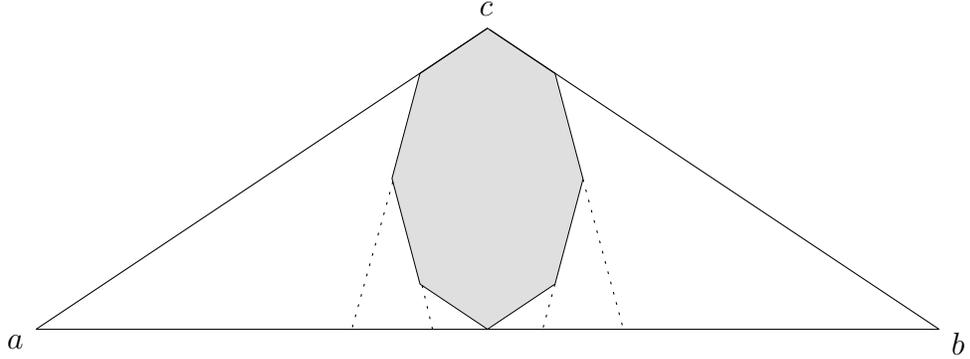
\lipefig{2}{Partitioning $G_1$ into triangles.}

\heading{Proof. }  Denote the vertices of $T$ by $a, b$ and $c$ so
that $U$ is a neighborhood of $ab$. We introduce a rectangular
coordinate system so that $ab$ lies on the $x$-axis and the
$y$-coordinate of $c$ is positive. The affine copies of $P$
used for the tiling will be referred to as {\it tiles}. We
construct the tiling by an algorithm describing the
successive tiles and their respective proper places. Let $G_i$
denote the closure of the untiled part of $T$ at the $i$-th
stage of the construction $(i=0,1,2,\ldots)$, at which point
$i$ tiles have been put in place. Obviously, at the
beginning, the number of tiles placed is $0$ and $G_0=T$. We
define the first tile, $P_1$, to be an affine copy of $P$ 
properly inscribed in $T$ and we partition $G_1$ into a
collection of triangles
${\cal T}_1=\left\{T_1^1,T_1^2,\ldots,T_1^{2n-2}\right\}$
each of which has its base on $ab$ and top vertex at some
vertex of the first tile (see Figure 2).

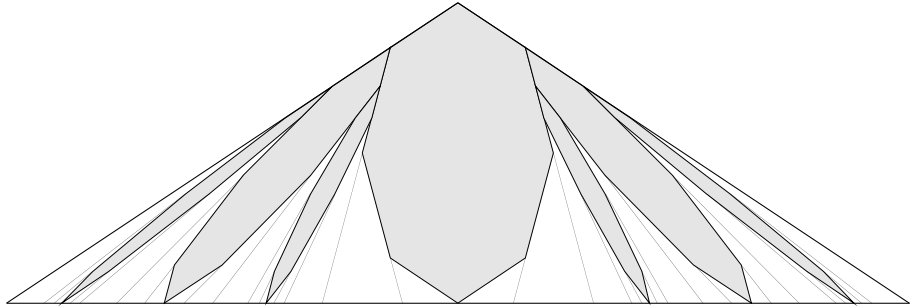
\lipefig{3}{Tiling a triangle minus its base's neighborhood.}

The formula for designing and placing the next ({\it i.e.},
the $(i+1)$-st) tile in $T$ is:

{\em Next Tile. } Among all triangles of ${\cal T}_i$ choose a
tallest one, {\it i.e.} one whose top vertex $v$ has a
maximum $y$-coordinate and call it $T_i^{\rm max}$. Let
$F$ be an affine transformation sending $T$ onto
$T_i^{\rm max}$ with $F(c)=v$ and define $P_{i+1}$ to be
$F(P_1)$, which is an affine copy of $P$ properly inscribed
in $T_i^{\rm max}$. Then define the partition ${\cal T}_{i+1}$
of $G_{i+1}$ by replacing $T_i^{\rm max}$ with the images of
the triangles in ${\cal T}_1$ under $F$. (Fig. 3 shows the
tiling stage at $i=7$).

Let now $y_i$ be the maximum of the $y$-coordinates of points
in the closure of $G_i$. Of course, $y_i$ occurs at one of
the vertices of $G_i$, thus at the top vertex of one of the
triangles of ${\cal T}_i$. Obviously, $y_i>0$, and
$y_{i+1}\le y_i$. Let $M_i$ be the line $y=y_i$ and let $m_i$
be the line $y=\mu y_i$, where $0<\mu<1$ is the constant
described in (2). As we place the successive tile at
the $(i+1)$-st stage of the construction, the top vertex of
the tile eliminates one vertex of $G_i$ lying on $M_i$, and,
by (2), every non-top vertex of this tile lies below the line
$m_i$. Thus, between the lines $m_i$ and $M_i$, the set of
vertices of $G_{i+1}$ is obtained from the set of vertices of
$G_i$ by deleting one element. Therefore there exists an
integer $k$ such that no vertex of $G_{i+k}$ lies above the
line $m_i$. It follows that for every $i$ there exists a $k$
such that $y_{i+k}\le \mu y_i$, and, consequently,
$\displaystyle\lim_{i\to\infty}y_i=0$. This implies that
there is an integer $k_0$ such that all vertices of $G_{k_0}$
lie in the neighborhood $U$ of $ab$, and the proof of the
lemma is complete.


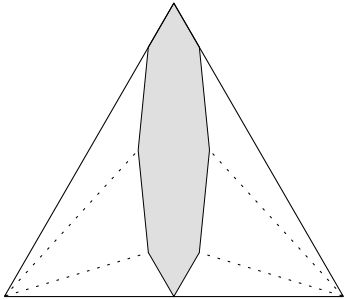
\lipefig{4}{Triangulating the complement of the initial tile.}

We now construct a periodic packing of the plane with
ellipses. Begin with the familiar regular tiling of the plane
with congruent equilateral triangles. Place in one of the
triangles a properly inscribed affine copy of $P$, and call
it the {\it initial tile}. Partition the remaining portion of
the triangle into $2n-2$ smaller triangles as shown in Fig.
4. The ellipse inscribed in the initial tile, when
homothetically $\lambda$-enlarged, covers a neighborhood of
the tile, hence it covers a neighborhood of one edge of each
of the smaller triangles. Tile each of the smaller triangles
minus a neighborhood of the edge already covered, in the
manner described in the proof of the Lemma (see Fig. 3).
Finally, extend this pattern to all triangles of the regular
tiling so that the initial tiles in each of them are
translates of each other (see Figure 5).

\lipefig{5}{A periodic arrangement.}

\heading{Remark (Istv\'an Talata). } Instead of ellipses, affine images
of an arbitrary centrally symmetric convex domain can be used
in the above construction, which would require minor
modifications only. In other words, the B\'ar\'any number of
the class of bounded-diameter affine images of a plane centrally symmetric
convex domain is zero.

\section{Packings with discs}
Throughout this section, we have $\eps=\epsval$.
Suppose for contradiction that there exists a packing $\CC$ with
discs of radius at most 1 such that its $(1+\eps)$-enlargement 
covers a square $S$ with side length $4$. Let us say that a disc $C\in\CC$
{\em bites\/} into a set $X\subseteq \R^2$ if $C^\eps\cap X\neq\emptyset$.
By induction, we are going to
construct a sequence of
compact sets $S=R_1\supset R_2\supset R_3\ldots$
such that for each $n=1,2,\ldots$, no disc of $\CC$ of radius
greater than $r_n$ bites
into $R_n$, where $(r_n)_{n=1}^\infty$ is a decreasing sequence of real
numbers tending to 0. Taking a point $x\in\bigcap_{n=1}^\infty R_n$
leads to a contradiction, since such an $x$ cannot be covered by
any $C^\eps$ with $C\in\CC$.

Each of the regions $R_n$ will be of one of two types,
called the {\em square type\/} and the {\em crescent type\/}.
We now describe the shape and the inductive hypothesis for these
two types of regions.

A region $R_n$  of the square type is a square of side $4r_n$,
and we assume that no disc $C\in \CC$ of radius greater than
$r_n$ bites into $R_n$. As a basis of the induction, we
choose $r_1=1$ and we let $R_1=S$.

A region $R_n$ of the crescent type is defined using some disc
$C_n\in \CC$, and $r_n$ is the radius of this $C_n$
(see Fig.~\ref{f:cr}).
\lipefig{cr}{A crescent-type region $R_n$.}
We fix suitable constants\footnote{The choice of the constants 
in the proof is somewhat arbitrary. The goal, rather than trying
to get the best value of $\eps$, was to select them in such a way
that realistic pictures can be drawn.}
 $\alpha=\alphaval$ and
$\beta=\betaval$. Let $c$ denote the center of $C_n$;
then $R_n$ is the intersection of an angle $\alpha$
with apex at $c$ with the $\beta$-ring of $C_n$.
We say that $R_n$ is a {\em crescent of\/} $C_n$,
and we call the halfline originating at $c$
and dividing $R_n$ into two equal parts
the {\em axis of $R_n$\/}.

We describe how $R_{n+1}$ is constructed from $R_n$.
First, we treat the simpler case when $R_n$ is of the square
type. Let $D$ be the disc of radius $\frac 32 r_n$ centered at the
center of the square $R_n$ (Fig.~\ref{f:sq}).
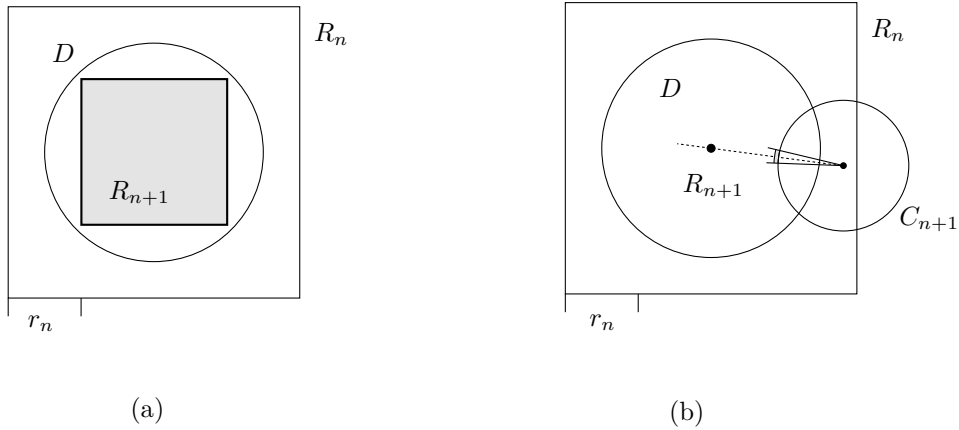
\lipefig{sq}{The inductive step for a square-type region $R_n$.}
Choose $C_{n+1}\in \CC$ as the disc of the largest radius
that bites into $D$. If the radius of $C_{n+1}$ is at most
$r_n/2$, set $r_{n+1}=r_n/2$ and
pick the region $R_{n+1}$ as a square of side $4r_{n+1}$
inside the disc $D$, as in Fig.~\ref{f:sq}(a).
Otherwise, let $r_{n+1}$ be the radius of $C_{n+1}$.
In this case, we pick $R_{n+1}$ as a crescent
of $C_{n+1}$. The axis of $R_{n+1}$ is
the halfline originating at the center of $C_{n+1}$
and passing through the center of $D$;
see Fig.~\ref{f:sq}(b) (if these centers happen to
coincide then pick an arbitrary  direction of the axis).
This finishes the definition of $R_{n+1}$. Easy geometric
considerations, whose details we omit, show that thus constructed
$R_{n+1}$ satisfies the inductive hypothesis (i.e. no disc of radius
larger than $r_{n+1}$ bites into $R_{n+1}$).

It remains to discuss the inductive step from $R_n$ to $R_{n+1}$
for an $R_n$ of the crescent type. 
For a simpler notation, we will measure distances in the units
of $r_n$ from now on, that is, we may assume $r_n=1$.
In this case, we let $I$
denote the intersection of $R_n$ with the $\Ifrac$-ring of
$C_n$ (see Fig.~\ref{f:i}(a)\,). 
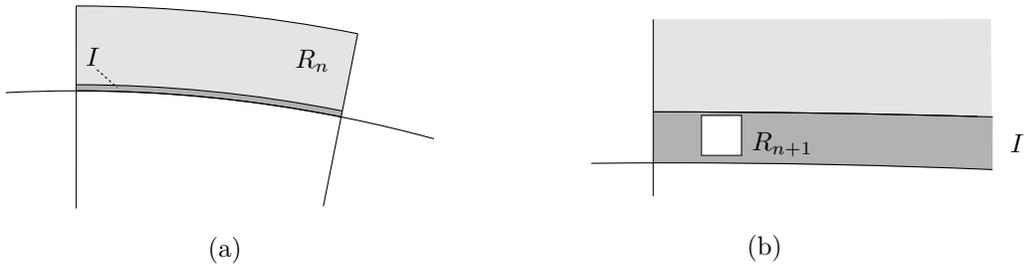
\lipefig{i}{The region $I$ (a), and the case of a very small $r$ (b).}
Let $ C$ be the largest disc of $\CC$
distinct from $C_n$
biting into $I$, and let $ r$ be the radius of $ C$.
Here we distinguish three cases: $ r\leq \smallrt $,
$\smallrt  <  r \leq \bigrt $, and $\bigrt < r \leq 1$.

\subh{The case $ r\leq\smallrt $. }
Here we set $r_{n+1}=\smallrt$ and we choose $R_{n+1}$ as a square
of side $\newsqs$ within $I$ so that $C_n$
doesn't bite into it, as in Fig.~\ref{f:i}(b). This is a valid
region of the square type.

\subh{The case $\smallrt < r \leq \bigrt $. }
Let $c_n$ denote the center of $C_n$, and let $c$
be the center of $C$. We set $C_{n+1}= C$,
we define $r_{n+1}$ as the radius of $C_{n+1}$, and we choose
$R_{n+1}$ as a crescent of $C_{n+1}$
as follows (Fig.~\ref{f:c2}).
\lipefig{c2}{The case $\smallrt < r\leq \bigrt$. }
The angle of the axis of $R_{n+1}$ and of the halfline
$cc_n$ is $\frac52 \alpha$ and $R_{n+1}$ lies on the side
of the segment $c_n c$ closer to the axis of $R_n$.
To verify the inductive hypothesis for $R_{n+1}$,
we first need to show that $R_{n+1}\subseteq R_n$.
That is, we need to check that the point
$v$ in Fig.~\ref{f:c2} has distance at most $1+\beta$
from $c_n$, that $u$ has distance at least $1$
from $c_n$, and that the point $w$
cannot go beyond the side boundary of $R_n$.
Let us check these conditions computationally.
As for the first inequality,
 $|vc_n|\leq 1+\beta$, we use the cosine theorem
for the triangle $c_ncv$:
$$
|c_nv|^2=|c_n c|^2+| cv|^2-2|c_n c|.| cv|\cos 3 \alpha =
(|c_nc|-|cv|)^2+2|c_nc|\cdot|cv|(1-\cos 3 \alpha)\le
$$
$$
\le (1+\frac\beta{16}+\frac\eps8)^2+2(1+\frac\beta{16}+(1+\eps)\frac18)
\cdot\frac18(1-\cos 3 \alpha)=1.055...
$$
Thus, $|c_n v| < 1.03 <1+\beta=1.0625$.

Similarly,
$$
|c_nu|^2=(|c_nc|-|cu|)^2+2|c_nc|\cdot|cu|(1-\cos 2\alpha)>
(1-\beta r)^2+2\cdot1\cdot r\cdot (1-\cos2\alpha)>
$$
$$
>1+2r(1-\cos2\alpha-\beta)> 1+\frac\beta{40}(1-\cos2\alpha-\beta)
=1.000021...
$$

To show that $R_{n+1}\subseteq R_n$, it remains to verify
that the point $w$
cannot go beyond the side boundary of $R_n$. Suppose to a contrary
that the angle $cc_nw$ is more than $\alpha\over2$. Then the distance
of the point $w$ to the line $cc_n$ is at least
$r_n\sin{\alpha\over2}=\sin{\alpha\over2}$. On the other hand,
the same distance equals $(1+\beta)r\sin 3\alpha$. It follows that
$\sin{\alpha\over2}\le(1+\beta)r\sin 3\alpha
\le(1+\beta){1\over8}\sin 3\alpha<{1\over6}\sin 3\alpha$,
which contradicts the concavity of the sine function in the interval
$[0,\pi]$.

To verify the induction hypothesis for $R_{n+1}$, it remains to show
that  no disc $C'$ of $\CC$ with radius in the
interval $(r,1]$ may bite into $R_{n+1}$.
Since $r\geq \smallrt$ is not too small, any such $C'$ biting into
$R_{n+1}$ would have to intersect $C_n$ or $C$. $C_n$ itself
doesn't bite into $R_{n+1}$, since $|c_nu|^2>1.000021$ and thus
$|c_nu|>1+\eps$.

\subh{The case $\bigrt<r\leq 1$. } Here we have the relatively
large disc $C$ biting into the region $I$ (Fig.~\ref{f:c3}).
\lipefig{c3}{The case $\bigrt<r\leq 1$.}
Consider the circular arc $a$ bounding the region $I$ from the outer
side, and let $a_0$ be the portion of this arc contained in the
the disc $C^\eps$. Calculation shows that even
if $C$ has the largest possible radius 1 and touches $C_n$
in the middle of the region $R_n$, one of the portions
of $a\setminus a_0$ has angular length at least $\newarc$
(this extreme case is shown in Fig.~\ref{f:c3}).
We thus select a portion $I'$ of the region $I$ of angular length
$\newarc$, avoiding $a_0$ but adjacent to it.

What is the largest possible radius of a disc
$C'\neq C_n$ of $\CC$ that may bite into $I'$?
The possible radius is largest when $r$ is smallest, that
is, $r=\bigrt$.  Fig.~\ref{f:c3a} shows how to  
upper-bound  the radius of $C'$; the radius of the disc $D$
drawn there, which is well below $\maxrprime$,
 is an upper bound for the radius of $C'$.
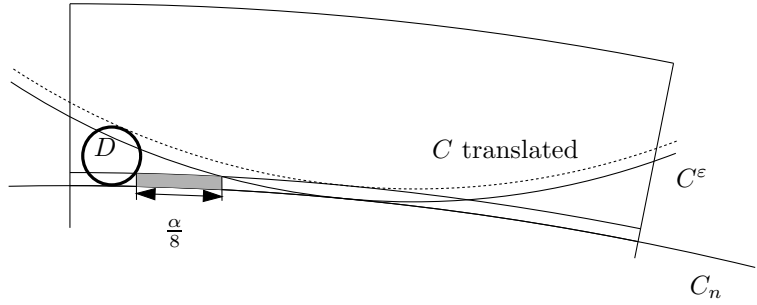
\lipefig{c3a}{Estimating the radius of $C'$.}
 Now we repeat the
considerations made above with the region $I'$
instead of $I$, that is, we choose the largest disc $C'\neq C_n$
biting into $I'$, we let $r'$ be its radius, and discuss
the cases depending on the range of $r'$. The first
two cases ($r'\leq \smallrt$ and $\smallrt < r' \leq \bigrt$) work
in the same way as above; the only small change is
that $I'$ is shorter than $I$, so one has to check
that there's always enough room to accommodate the
corner $w$ of the region $R_{n+1}$ (as in Fig.~\ref{f:c2})
in the region $R_n$. But this works because $r'$
is small enough.
And, because of the restriction $r'\leq \maxrprime$,
the third case discussed for the region $I$ cannot occur
for $I'$.
Theorem~\ref{t:discs} is proved.


\end{document}